\documentclass[final]{amsart}
\usepackage{amssymb,amsmath, color}

\usepackage{amssymb}

\newtheorem{Theorem}{Theorem}[section]
\newtheorem{lemma}{Lemma}[section]
\newtheorem{corollaire}{Corollary}[section]

\theoremstyle{remark}
\newtheorem{remark}{Remark}[section]

\newcommand{\R}{{\mathbb R}}

\newcommand{\argmin}{{\rm argmin}\kern 0.12em}

\newcommand{\X}{\mathcal X}

\newcommand{\Z}{\mathcal Z}
\newcommand{\Hb}{\mathcal H}

\newcommand{\Xiu}{\mathcal X_1}
\newcommand{\Xid}{\mathcal X_2}

\usepackage{ulem}

\begin{document}

\title{Asymptotic behavior of gradient-like dynamical systems involving inertia and multiscale aspects
  }

\author{Hedy Attouch}

\author{Marc-Olivier Czarnecki}

\address{Institut Montpelli\'erain Alexander Grothendieck, UMR 5149 CNRS, Universit\'e Montpellier 2, place Eug\`ene Bataillon,
34095 Montpellier cedex 5, France}
\email{hedy.attouch@univ-montp2.fr, 
marc-olivier.czarnecki@univ-montp2.fr}

 \date{\today}     
 \thanks{H. Attouch: with the support of  ECOS  grant C13E03. Effort sponsored by the Air Force Office of Scientific Research, Air Force Material Command, USAF, under grant number FA9550-14-1-0056.}
\maketitle

\vspace{0.3cm}

\paragraph{\textbf{Abstract}} In a Hilbert space $\mathcal H$, we study the asymptotic behaviour, as  time variable $t$ goes to $+\infty$, of  nonautonomous gradient-like dynamical systems involving inertia and multiscale features.
 Given 
% $\mathcal H$ a general Hilbert space}, 
 $\Phi: \mathcal H \rightarrow \R$ and $\Psi: \mathcal H \rightarrow \R$  two convex differentiable functions, 
 $\gamma$  a positive damping parameter, and $\epsilon (t)$ a function of $t$ which tends to zero as $t$ goes to $+\infty$,
we consider the  second-order differential equation
$$
\ddot{x}(t)   + \gamma  \dot{x}(t) +    \nabla \Phi (x(t)) +   \epsilon (t) \nabla \Psi  (x(t)) = 0. 
$$
This system  models the emergence of various  collective behaviors in game theory, as well as the asymptotic control of  coupled nonlinear oscillators. 
Assuming that $\epsilon(t)$ tends to zero moderately slowly as $t$ goes to infinity, we show that the trajectories  converge weakly in $\mathcal H$. The limiting equilibria are solutions of the hierarchical minimization problem which consists in minimizing
$\Psi$ over the set  $C$ of minimizers of $\Phi$.  As  key assumptions, we 
suppose that $  \int_{0}^{+\infty}\epsilon (t) dt = + \infty  $ and that, for every $p$ belonging to a convex cone $\mathcal C$ depending on the data $\Phi$ and $\Psi$
$$  \int_{0}^{+\infty}  \left[\Phi^* \left(\epsilon (t)p\right) -
\sigma_C \left(\epsilon (t)p\right)\right]dt < + \infty
$$
where $\Phi^*$ is the Fenchel conjugate of $\Phi$, and     $\sigma_C $ is the support function of $C$.
An application is given to coupled oscillators.

\vspace{0.3cm}


\vspace{0.3cm}


\markboth{H. ATTOUCH, M.-O. CZARNECKI}
  {MULTISCALED GRADIENT DYNAMICS}

\indent

\newpage

\section{Introduction}

$\mathcal H$ is a real Hilbert space, we write $ \|x\|^2 = \left\langle x , x\right\rangle$ for $x\in \mathcal H$. For any differentiable function $F: \mathcal H \rightarrow \R $, its gradient is denoted by $\nabla F$. Thus $F'(x)(y) =  \left\langle \nabla F(x) , y\right\rangle$. The first order (respectively second order) derivative at time $t$ of a function $x(\cdot): [0,+\infty [ \to \mathcal H$ is denoted by $\dot{x}(t)$ (respectively $\ddot{x}(t)$). Throughout the paper $\gamma$ is a fixed positive parameter (viscous damping coefficient).

\subsection{Problem statement} 

Henceforth, we make the following standing assumptions $(\mathcal H_0)$ on data $\Phi$, $\Psi$ and $\epsilon(\cdot)$, that will be needed throughout the paper.  We write $S=\mbox{\argmin}_C \Psi$ with $C$ = $\mbox{\argmin}\Phi$, or, equivalently,  $S=\mbox{\argmin}\{\Psi|\mbox{\argmin}\Phi\}$.

\medskip

$(\mathcal H_0)\quad\left\{\rule{0em}{3em}\right.$
\parbox{45em}{
  $\Phi: \mathcal H \rightarrow \R^{\phantom{+}} $ is   convex,  $\nabla \Phi$ is Lipschitz continuous on bounded sets of $\mathcal H$, 
$C$ = $\argmin \Phi = \Phi ^{-1}(0) \neq \emptyset$;\\
  $\Psi: \mathcal H \rightarrow \R $  is  convex,  $\nabla \Psi$ is Lipschitz continuous on bounded sets of $\mathcal H$, $\Psi$ is bounded from below;\\
  $S=\mbox{\argmin}\{\Psi|\mbox{\argmin}\Phi\} \neq \emptyset$;\\
  $\epsilon(\cdot): [0,+\infty[ \to ]0,+\infty [ \mbox{ }  \mbox{is a}\mbox{ }   \mbox{ nonincreasing function,}  \mbox{ of class } C^1 , \    \epsilon(t) \mbox{ }\mbox{and} \mbox{ }\dot{\epsilon}(t) \mbox{ }\mbox{tend to zero   as} \mbox{ } t\to +\infty$. 
}

\medskip

We study the asymptotic behavior ($t \rightarrow + \infty  $)  of the trajectories of the nonautonomous Multiscaled Inertial Gradient-like  system ((MIG) for  short)
\begin{equation}\label{MIG}
\tag{\rm MIG}   \qquad     \ddot{x}(t)   + \gamma  \dot{x}(t) +    \nabla \Phi (x(t)) +   \epsilon (t) \nabla \Psi  (x(t)) = 0. 
\end{equation}

\vspace{0.1cm}

 Let us observe that when  $\Psi \equiv 0$, or $\epsilon =0$, the above dynamic reduces to
\begin{equation}\label{HBF}
\tag{\rm HBF}   \qquad     \ddot{x}(t)   + \gamma  \dot{x}(t) +    \nabla \Phi (x(t))  = 0. 
\end{equation}
The Heavy Ball with Friction dynamical system, \eqref{HBF} for short, plays an important role in mechanics, control theory, and optimization, see \cite{Al}, \cite{ACR}, \cite{AGR}, \cite{HJ}, \cite{Polyak} for a general presentation.
%The (HBF) system plays an important role in mechanics, control theory, and %optimization, see \cite{AGR}.  It is often called the Heavy Ball with Friction %dynamical system, hence the terminology (HBF).
Because of the viscosity friction term $\gamma  \dot{x}$, it is a dissipative dynamical system (the global mechanical energy is decreasing), which gives it remarkable properties of optimization. Indeed, just assuming that $\Phi$ is  convex and $C= \argmin \Phi \neq  \emptyset$,       each trajectory of system \eqref{HBF}  converges weakly in $\mathcal H$, with a limit belonging  to  the closed convex set $C$, see Alvarez \cite{Al}.

In many situations, there is a continuum of equilibria (equivalently, $C$ is not reduced to a singleton), as it occurs when considering  a potential $\Phi$ which is not strictly convex. A  situation frequently encountered results from the application of a least squares method to an under determined linear problem. Another important situation that we address in this paper concerns the weakly coupled systems.
In this paper, we address the question: how to control \eqref{HBF} to obtain asymptotically an equilibrium having desirable properties?
A somewhat related issue is: why specific equilibria are observed in real life, despite the fact that no force is present who can explain?

We give an answer to these questions  by introducing in \eqref{HBF} an asymptotically vanishing term $\epsilon (t) \nabla \Psi  (x(t))$, so obtaining \eqref{MIG}.  The crucial point  is to take a control variable  $\epsilon (t)$ which  tends to zero in  a \emph{moderate} way, that is,  not too fast
 $$\int_0^{+\infty} \epsilon (t) dt = + \infty, $$
  and not too slow: we suppose that the following inequality holds for every $p$ belonging to a cone, whose definition involves $\Phi$ and $\Psi$, 
$$
\int_{0}^{+\infty}  \left[\Phi^* \left(\epsilon (t) p\right) - \sigma_C\left(\epsilon (t) p\right)\right]dt 
< + \infty,
$$ 
where $\Phi^*$ is the Fenchel conjugate of $\Phi$, and $\sigma_C$ is the support function of $C$.
In this case,  we will prove the following asymptotical hierarchical selection property: each trajectory of system  \eqref{MIG} converges weakly, with a limit belonging  to  $\argmin_C \Psi$, i.e., the limit does minimize $\Psi$ over the set of minimizers of $\Phi$:
\begin{equation}\label{MAG3}
x(t)\rightharpoonup x_{\infty} \in \argmin_C \Psi  \mbox{ as }  t\rightarrow +\infty.
\end{equation}
The case $\Psi  (x) =  \|x\|^2$, which corresponds to Tikhonov regularization, has been examined by Attouch and Czarnecki in \cite{AtCz1}. In this case, each trajectory of  \eqref{MIG} converges strongly  to the unique element of $C$ with minimal norm. A natural extension of this result to the case of a uniformly convex function $\Psi$  was obtained
by Cabot in \cite{Cabot-inertiel}.

For many applications, it is important to go beyond  the strongly and uniformly convex cases. A typical situation, that we illustrate in Section~\ref{coupled_section}, 
 is the  modeling of the (weak) coupling between two oscillators $x(\cdot)$ and $y(\cdot)$, where we take  $\Psi (x,y) =\|Ax - By\|^2$, $A$ and $B$ linear  continuous operators.
We also give an application of this model  to the Nash equilibration for potential games, and to coupled wave equations.

% Returning to system (\ref {MIG}), the same property is satisfied by the solution trajectories of \eqref{MIG}, under the assumption that $\epsilon (\cdot)$ tends to zero sufficiently fast, i.e., 
% $\int_0^{+\infty} \epsilon (t) dt < + \infty$. 

\noindent The asymptotic analysis of the  multiscaled  first-order differential inclusion 
\begin{equation*}\label{MAG}
   \qquad     \dot{x}(t)   + \partial \Phi (x(t)) +    \epsilon(t) \partial  \Psi (x(t))  \ni 0, 
\end{equation*}
involving general convex potentials $\Phi$ and $\Psi$ has been studied by  Attouch and Czarnecki in \cite{AtCz2}. As we shall see, several asymptotic convergence properties can be carried from the first order to the second order differential system.
In particular, the moderate decrease condition on $\epsilon(\cdot)$ which allows to obtain an asymptotic hierarchical minimization property, and which involves the Fenchel conjugate of $\Phi$ is 
quite similar. Note  that, by contrast with the first-order case where general convex lower semicontinuous potentials are considered, our analysis of the inertial second-order differential system deals with differentiable potentials. As a general rule, introducing nonsmooth potentials in inertial systems leads us to consider trajectories whose acceleration 
is a vectorial measure, and with possible shocks, an interesting topic whose study goes beyond the scope of this article, see \cite{ACR} for the (HBF) system.

 The local Lipschitz continuity assumption  on the gradient of the potentials $\Phi$ and $\Psi$ ensures the existence of strong solutions to (MIG), but this assumption is not used in the  proof of convergence, suggesting some extensions to PDE's. An illustrative example is given in Section \ref{coupledwavesection} concerning coupled wave equations.  In this case, the potential $\Phi$ (Dirichlet energy)  is only lsc. on $\mathcal H = L^2 (\Omega)$. Still, the existence can be proved by direct energy methods, and we show that the result of convergence holds.

\subsection{Moderate decrease condition on $\epsilon(\cdot)$.}

To prove the convergence of trajectories of  \eqref{MIG} system, and show the asymptotic hierarchical minimization property (\ref{MAG3}), we  need to make assumptions about how  $\epsilon(\cdot)$ 
tends to zero as $t$ goes to infinity. These   conditions
express that  $\epsilon(\cdot)$   is not too small $(\mathcal H_1) $, and not too large $(\mathcal H_2)$:
\begin{equation}
\tag{$\mathcal H_1$}	    \int_0^{+\infty} \epsilon (t) dt = + \infty; 
\end{equation}
\begin{equation}
\tag{$\mathcal H_2$}	  \forall p\in \mathcal C \mbox{ } \mbox{ }\mbox{ } \int_{0}^{+\infty}  \left[\Phi^* \left(\epsilon (t) p\right) - \sigma_C\left(\epsilon (t) p\right)\right]dt 
< + \infty,
\end{equation}
with
\begin{equation}\label{def-cone}
 \mathcal C  =  \left\{ \lambda p:  \ \lambda \geq 0,  \exists z\in C, \  p = - \nabla \Psi (z) ,  \    \nabla \Psi (z)+ N_C (z) \ni 0  \     \right\}.
\end{equation}
\noindent In $(\mathcal H_2)$ we use classical concepts and notations of convex analysis: $\Phi^*$ is the Fenchel conjugate of $\Phi$, defined for every $y\in  \mathcal H $ by 
$$
\Phi^* (y) = \sup  \left\{  \langle  y ,  x \rangle - \Phi (x) |{x\in   \mathcal H}\right\},
$$
\noindent $ \sigma_C $ is the support function of $C= \argmin \Phi$ , defined for every $y\in  \mathcal H $, by 
$$\sigma_C (y) = \sup  \left\{ \langle  y ,  x \rangle|{x\in C}\right\},$$
\noindent and $ N_C (x)$ is the (outwards) normal cone to $C$ at $x$.\\
The cone $\mathcal C$ can be equivalently  defined as
$$
\mathcal C  =  \left\{ \lambda p:  \ \lambda \geq 0,  \exists z\in \argmin_C \Psi , \  p = - \nabla \Psi (z)   \     \right\}.
$$
One can verify that, under some general qualification assumption, $\mathcal C$ is a closed convex cone. Indeed
$$
\mathcal C  =  \left\{ \lambda p:   \lambda \geq 0, p\in \argmin_{\mathcal H}\left(\sigma_C(\cdot) + \Psi^*(- \cdot)\right)\right\}.
$$

\vspace{0.2cm}
\textit{Analysis of  condition} $(\mathcal H_1)$: the slow parametrization  condition 
$ \int_0^{+\infty} \epsilon (t) dt = + \infty$ expresses that the control variable $\epsilon (\cdot)$ converges slowly to zero, when time $t$ goes to infinity. It plays a central role in the  asymptotic hierarchical property. This fact has been well established in a series of papers concerning  first or second order systems of this type, see \cite{ACom}, \cite{AtCz1}, \cite{AtCz2}, \cite{Cabot}, \cite{Cabot-inertiel}.
\\
In contrast, the fast parametrization case $ \int_0^{+\infty} \epsilon (t) dt < + \infty$,  behaves exactly as if $\epsilon \equiv 0$. In this case, there is no asymptotic hierarchical property associated to the control variable.

\vspace{0.2cm}
\textit{Analysis of  condition} $(\mathcal H_2)$:
\vspace{0.1cm}
\begin{itemize}
\item a) Note that $\Phi$ enters in \eqref{MIG} only via its subdifferential. Thus it is not a restriction to assume $\min_{\mathcal H} \Phi =0$. 
For a function $\Phi$ whose minimum is not equal to zero, one should replace in $(\mathcal H_2)$ and in the corresponding statements $\Phi$ by $\Phi - \min_{\mathcal H} \Phi.$

\noindent From $\Phi \leq {\delta}_C$ we get  $\Phi^* \geq \left(\delta_C\right)^* =  \sigma_C$,   and $ \Phi^*  - \sigma_C  \geq 0.$

\noindent We have used that $ \sigma_C $ is equal to the Fenchel conjugate of ${\delta}_C$, where ${\delta}_C$ is the indicator function of $C$.

\noindent Hence,  $(\mathcal H_2)$ means that, for $p\in \mathcal C$,  the  nonnegative function 
$$ t \mapsto  \left[\Phi^*  \left(\epsilon (t) p\right) - \sigma_C\left(\epsilon (t)p\right)\right]$$ 
\noindent is integrable on $[0, +\infty[$. Practically, this condition forces $\epsilon(\cdot)$  to be small enough for large $t$.

\vspace{0.1cm}

\item b) As an illustration, consider the model situation: $\Phi (z) =  \frac{1}{2} \mbox{dist}^2 (z, C) =  \frac{1}{2} \|.\|^2  \mbox{ }+_e \mbox{ } \delta_C$, where $+_e$ denotes the epigraphical sum (also called inf-convolution). From  general properties of  Fenchel transform 

\begin{center}
$\Phi^* (z) =  \frac{1}{2} \|z\|^2  +  \sigma_C (z)$ and $\Phi^* (z) -   \sigma_C (z)  =  \frac{1}{2} \|z\|^2 $. 
\end{center}

\noindent Hence, in this situation

$$(\mathcal H_2) \Longleftrightarrow  \int_{0}^{+\infty} \epsilon(t)^2 dt < + \infty.$$

\noindent Conditions $(\mathcal H_1)$  and $(\mathcal H_2)$    are equivalent to $ \epsilon (\cdot) \in L^2 (0, +\infty) \setminus L^1 (0, +\infty)$.
It is satisfied for example by taking $\epsilon (t) = \frac{1}{1+ t^{\alpha}}$ \ with    \   $\frac{1}{2} < \alpha \leq 1$.

\end{itemize}

\subsection{Contents}
In Section \ref{section_sub}, we express our main result, which shows that, under conditions of moderate decrease of $ \epsilon (\cdot) $, each trajectory of \eqref{MIG} converges weakly to the solution of a hierarchical minimization  problem. In Section~\ref{section_strong}, we complement this result by showing a property of strong convergence when $ \Psi $ is strongly convex.
Then, in Section~\ref{multiscale},  we look at some multiscale aspects, and show the effect of changing the scaling in dynamic. In particular, we make the link between asymptotic vanishing viscosity and penalization. Finally, we illustrate our study with applications to coupled oscillators in Section~\ref{coupled_section}, and coupled wave equations in  Section~\ref{coupledwavesection}.

\section{Weak convergence results}\label{section_sub}

In this section, we consider the dynamical system
\begin{equation*}
\tag{\ref{MIG}} \qquad \ddot{x}(t)   + \gamma  \dot{x}(t) +    \nabla \Phi (x(t)) +   \epsilon (t) \nabla \Psi  (x(t)) = 0.
\end{equation*}
\noindent When $\Psi=0$, \eqref{MIG} boils down to the classical second order damped nonlinear oscillator
\begin{equation*}\tag{HBF}\ddot{x}(t)   + \gamma  \dot{x}(t) +    \nabla \Phi (x(t)) = 0.
\end{equation*}
In accordance with Alvarez theorem \cite{Al},  and Attouch-Czarnecki \cite {AtCz1},    in our main
result, which is stated below, we are going to show that, under moderate decrease condition on $\epsilon(\cdot)$,  each trajectory of \eqref{MIG} converges  weakly in $\Hb$ to a
minimizer of $\Phi$, which also minimizes $\Psi$ over all minima of
$\Phi$.

\begin{Theorem} \label{basicthm}
Let us suppose that standing assumptions $(\mathcal H_0)$ are satisfied:

$(\mathcal H_0)\  \left\{\rule{0em}{3em}\right.$
\parbox{45em}{
  $\Phi: \mathcal H \rightarrow \R^{\phantom{+}} $ is   convex,  $\nabla \Phi$ is Lipschitz continuous on bounded sets of $\mathcal H$, 
$C$ = $\argmin \Phi = \Phi ^{-1}(0) \neq \emptyset$;\\
  $\Psi: \mathcal H \rightarrow \R $  is  convex,  $\nabla \Psi$ is Lipschitz continuous on bounded sets of $\mathcal H$, $\Psi$ is bounded from below;\\
  $S=\mbox{\argmin}\{\Psi|\mbox{\argmin}\Phi\} \neq \emptyset$;\\
   $\epsilon(\cdot): [0,+\infty[ \to ]0,+\infty [ \mbox{ }  \mbox{is a}\mbox{ }   \mbox{ nonincreasing function,}  \mbox{ of class } C^1 , \    \epsilon(t) \mbox{ }\mbox{and} \mbox{ }\dot{\epsilon}(t) \mbox{ }\mbox{tend to zero   as} \mbox{ } t\to +\infty$.  
}

Let $x$ be a classical maximal solution of \eqref{MIG}. Then $x(\cdot)$ is defined on $[0, +\infty [$.
Let us assume moreover that the moderate decrease conditions on $\epsilon(\cdot)$ are satisfied:
\begin{align}
	& \tag{$\mathcal H_1$}     \int_0^{+\infty} \epsilon (t) dt = + \infty ;\\
	&\tag{$\mathcal H_2$}     \forall p\in \mathcal C  \mbox{ } \int_{0}^{+\infty}  \left[\Phi^* \left(\epsilon (t) p\right) - \sigma_C\left(\epsilon (t) p\right)\right]dt 
< + \infty, \\
&\mbox {with }  \mathcal C  =  \left\{ \lambda p:  \ \lambda \geq 0,  \exists z\in C, \  p = - \nabla \Psi (z) ,  \    \nabla \Psi (z)+ N_C (z) \ni 0  \     \right\};
\nonumber\\
 & \tag{$\mathcal H_3$}    0 \leq - \dot{\epsilon}(t) \leq k \epsilon^2 (t)  \  \mbox{ }\mbox{ }  \mbox{for some positive constant} \  k, \mbox{ }  \mbox{and large enough} \mbox{ } t.
\end{align}

Then
\begin{eqnarray*}
(i)&\mbox{ weak convergence }&\exists x_{\infty}\in S=\mbox{\argmin}\{
  \Psi|\mbox{\argmin}\Phi\},\ \  \mbox{such that  }   w-\lim_{t\to+\infty}
  x(t)=x_{\infty};\\
(ii)&\mbox{ minimizing properties } &\lim_{t\to+\infty} \Phi(x(t))=0    \ \    \mbox{and}  \  \      \lim_{t\to +\infty} \Psi(x(t)) = \min\Psi |_{\mbox{\argmin}\Phi}.
\end{eqnarray*}
\end{Theorem}

Taking $\Psi=0$ in  Theorem~\ref{basicthm}, we recover the
 Alvarez convergence result for the heavy ball with friction system.

\begin{corollaire} \label{Alvarezthm}{\cite[ Theorem 2.1]{Al}}
Let
$\Phi: \mathcal H \rightarrow \R$ be a convex  function whose gradient $\nabla \Phi$ is Lipschitz continuous on the bounded subsets of $\mathcal H$, and such that
$C$ = $\argmin\Phi \neq \emptyset$. Let   $\gamma >0$ be a positive damping parameter.   Then for any solution trajectory  $x(\cdot)$ of
\begin{equation*}
\tag{\ref{HBF}}   \ddot{x}(t)   + \gamma  \dot{x}(t) +    \nabla \Phi (x(t))  = 0,
\end{equation*}
 $x(t)$  converges weakly in $\Hb$ to a point in $\argmin\Phi$, as t goes to $+\infty$.
\end{corollaire}

\textbf{Proof of corollary \ref{Alvarezthm}.}  Take $\Psi\equiv 0$, and choose $\epsilon(\cdot)$ an arbitrary function satisfying the standing assumptions, $(\mathcal H_1$) and $(\mathcal H_3)$
 (for example take $\epsilon(t)= \frac{1}{1+t}$). We claim that $(\mathcal H_2)$ is automatically satisfied. Since $\Psi\equiv 0$, the cone $\mathcal C$ is reduced to the origin. Noticing that $ \Phi^* (0) = - \inf_{\mathcal H} \Phi = 0$, and that $\sigma_C (0) = 0$, we have 
$\int_{0}^{+\infty}  \left[\Phi^*  \left(\epsilon (t) 0\right) - \sigma_C\left(\epsilon (t)0\right)\right] dt= 0$. Thus $(\mathcal H_2)$ is satisfied,  and we can apply Theorem \ref{basicthm}.
As a conclusion, we obtain the weak convergence of $x(t)$  to a point in $\argmin\Phi$, as $t$ goes to $+\infty$ (since $\Psi$ is constant, there is no hierarchical minimization associated to it).
That's the conclusion of Alvarez theorem.

\begin{remark}
The counterexample of Baillon~\cite{Baillon} shows that one may not
have strong convergence for \eqref{HBF}. Of course, the same holds for \eqref{MIG}.
\end{remark}

\subsection*{Proof of Theorem~\ref{basicthm}}
As in Bruck~\cite{Bruck}, weak convergence is a consequence of
Opial's lemma, after showing the convergence of $ \|x(.) -z\|$ for
every $z\in S$, and that every weak cluster point of $x$ belongs to
$S$. The proof consists of three steps, each of these steps relying on a different Lyapunov function.

%\medskip

\subsubsection*{Step 1:  First energy estimates and global existence results}

First note that, for any Cauchy data $x(0) = x_0$, $\dot{x}(0) = \dot{x}_0$, with $x_0 \in  \mathcal H $, $\dot{x}_0\in  \mathcal H $,  local existence and uniqueness of the corresponding classical solution to \eqref{MIG} is a direct consequence of Cauchy-Lipschitz theorem, after reducing \eqref{MIG} to a first order system.

Let $x :  [0, T[ \to \mathcal H $ \  be a maximal solution of \eqref{MIG}. To prove that $T = + \infty$ we use a classical contradiction argument. Let us suppose that $T<+\infty$.  To obtain a contradiction, it is enough to show that $\lim_{t\to T} x(t)$  and 
$\lim_{t\to T} \dot{x}(t)$ exist. This is a consequence of the energy estimates which are described below, which show that $\|\dot{x}(t)\|$ and $\|\ddot{x}(t)\|$ are bounded on $[0, T[$.
%{\color{red} conclusion sur l'existence?}
In order to avoid repeating twice the same argument, we  establish these estimations directly on global solutions of (MIG), with  $T= + \infty$.

 Given $x(\cdot):  [0, +\infty[ \to \mathcal H $ a global classical solution of \eqref{MIG}, set
\begin{equation}\label{Energ1}
E_1(t) =   \frac{1}{2} \|\dot{x}(t)\|^2   +     \Phi(x(t))   +  \epsilon (t)\Psi (x(t)).
\end{equation}
We are going to show that $E_1$ is a Lyapunov-like function, which will allow us to derive several asymptotic properties of trajectories of \eqref{MIG}. 
\begin{lemma}\label{energy1}  Let us just assume that standing assumptions $(\mathcal H_0)$ hold. Then, for any trajectory $x(\cdot):  [0, +\infty[ \to \mathcal H $  of \eqref{MIG}, the following properties hold:
\begin{itemize}
\item (i) $t \mapsto E_1(t) - c \epsilon(t)$ is a decreasing function, with $ c = \inf_ {\mathcal H} \Psi$;
\item (ii) $\dot{x} \in L^2( [0, +\infty[; \mathcal H) \cap L^{\infty} ( [0, +\infty[; \mathcal H) $.
\end{itemize} 
Assuming moreover that $x(\cdot)$ is bounded, then  
\begin{itemize}
\item (iii) $\lim_{t \to  +\infty} \dot{x}(t) = \lim_{t \to  +\infty} \ddot{x}(t)=0$;
\item (iv) $\lim_{t \to  +\infty} \nabla \Phi(x(t))= 0,$ \  and hence \   $\lim_{t \to  +\infty} \Phi(x(t))= \inf_{\mathcal H} \Phi$.
\end{itemize} 
\end{lemma}
\textbf{Proof of Lemma \ref{energy1}}.
Let us compute the time derivative of $E_1(\cdot)$. Using the classical derivation chain rule and  equation \eqref{MIG}, we obtain
\begin{eqnarray*}
\dot{E_1}(t)&=&\langle \dot{x}(t), \ddot{x}(t)\rangle   +   \langle \nabla \Phi (x(t)), \dot{x}(t)\rangle  +  \epsilon (t)\langle \nabla \Psi (x(t)), \dot{x}(t)\rangle  
+ \dot{\epsilon}(t) \Psi (x(t)) \\
&=&\langle \dot{x}(t), \ddot{x}(t)   +  \nabla \Phi (x(t)) +  \epsilon (t) \nabla \Psi (x(t)) \rangle   + \dot{\epsilon}(t) \Psi (x(t)) \\
&=& - \gamma \|\dot{x}(t)\|^2  + \dot{\epsilon}(t) \Psi (x(t)).
\end{eqnarray*}
Using that $\epsilon(\cdot)$ is a decreasing function (i.e., $\dot{\epsilon}(t) \leq 0$), and that $\Psi$ is bounded from below, we deduce that
\begin{equation*}\label{Energ12}
\dot{E_1}(t)     +\gamma \|\dot{x}(t)\|^2   \leq   c \dot{\epsilon}(t),
\end{equation*}
with $c =\inf_ {\mathcal H} \Psi $, that's item $(i)$.
After integration of this inequality with respect to $t$, we get 
\begin{eqnarray}\label{Energ131}
E_1(t) - E_1(0)  +\gamma \int_0^t \|\dot{x}(s)\|^2 ds     &\leq&  c ( \epsilon (t) -  \epsilon (0) )      \\
&\leq &  |c| \epsilon (0). \nonumber
\end{eqnarray}
As a consequence
\begin{equation*}\label{Energ13}
\sup _{t\geq 0} E_1(t) < + \infty.
\end{equation*}
By definition (\ref{Energ1}) of $E_1$, we infer the existence of some constant $C_1$ such that for all $t\geq 0$
\begin{equation*}\label{Energ14}
\frac{1}{2} \|\dot{x}(t)\|^2   +     \Phi(x(t))   +  \epsilon (t)\Psi (x(t)) \leq C_1.
\end{equation*}
Since $\Phi$ and  $\Psi$ are bounded from below and $\epsilon(\cdot)$ is bounded, it follows that
\begin{equation}\label{Energ15}
\sup _{t\geq 0} \|\dot{x}(t)\| < + \infty.
\end{equation}
Noticing that $E_1(\cdot)$ is bounded from below, inequality (\ref{Energ131}) also yields
\begin{equation}\label{Energ16}
\int_0^{+\infty} \|\dot{x}(t)\|^2 dt  < + \infty.
\end{equation}
Collecting (\ref{Energ15})   and    (\ref{Energ16}),  we  obtain  item  $(ii)$.\\
Let us now suppose that $x(\cdot)$ is bounded, i.e., 
\begin{equation*}\label{Energ17}
\sup _{t\geq 0} \|x(t)\| < + \infty.
\end{equation*}
By standing assumptions $(\mathcal H_0)$,  $\nabla \Phi$ and $\nabla \Psi$ are Lipschitz continuous on bounded sets of $\mathcal H$, and hence bounded on bounded sets of $\mathcal H$.
From equation \eqref{MIG}, also using the boundedness of $\epsilon(\cdot)$, and the boundedness of $\dot{x}(\cdot)$ (item $(ii)$), it follows that 
\begin{equation*}\label{Energ18}
\ddot{x}(\cdot)  \in L^{\infty} ([0, +\infty); \mathcal H).
\end{equation*}
Hence $\dot{x}(\cdot)$ is Lipschitz continuous on $]0, +\infty[$, a property which combined with $\dot{x} \in L^2 ( ]0, +\infty[; \mathcal H)$ classically implies
\begin{equation}\label{Energ19}
\lim_{t \to  +\infty} \dot{x}(t) = 0.
\end{equation}
Let us now prove that 
\begin{equation}\label{Energ20}
\lim_{t \to  +\infty} \ddot{x}(t) = 0.
\end{equation}
For the sake of simplicity, we assume that $\Phi$ and $\Psi$ are twice differentiable (otherwise, the argument can be justified by using finite difference quotients).
After derivation of equation \eqref{MIG} we obtain
\begin{equation*}\label{Energ21}
\dddot{x}(t)   + \gamma  \ddot{x}(t) +    \nabla^2 \Phi (x(t))\dot{x}(t) +  \epsilon (t) \nabla^2 \Psi (x(t))\dot{x}(t) +   \dot{\epsilon} (t) \nabla \Psi (x(t))= 0. 
\end{equation*}
The Lipschitz continuity on bounded sets of $\nabla \Phi$ and $\nabla \Psi$, together with (\ref{Energ19}) imply 
\begin{equation*}\label{Energ22}
\lim_{t \to  +\infty} \nabla^2 \Phi (x(t))\dot{x}(t)= \lim_{t \to  +\infty}   \nabla^2 \Psi (x(t))\dot{x}(t) = 0. 
\end{equation*}
Using too that $\lim_{t \to  +\infty}\epsilon(t) = \lim_{t \to  +\infty}\dot{\epsilon}(t)= 0$, we obtain 
\begin{equation}\label{Energ23}
\dddot{x}(t)   + \gamma  \ddot{x}(t) =  g(t) 
\end{equation}
for some function $g$ such that $ \lim_{t \to  +\infty} g(t) = 0$.
Integration of (\ref{Energ23})  yields 
\begin{equation*}\label{Energ24}
\ddot{x}(t)  = e^{-\gamma t} \ddot{x}(0) + e^{-\gamma t} \int_0^t e^{\gamma s} g(s) ds,
\end{equation*}
which, by $ \lim_{t \to  +\infty} g(t) = 0$, and an elementary argument, gives (\ref{Energ20}).\\
Let us now return to \eqref{MIG} equation to obtain
\begin{equation}\label{Energ25}
\lim_{t \to  +\infty} \nabla \Phi (x(t)) = 0. 
\end{equation}
By using a  standard convexity argument we are going to deduce that
\begin{equation}\label{Energ26}
\lim_{t \to  +\infty} \Phi (x(t)) = \inf_{\mathcal H} \Phi. 
\end{equation}
To that end, let us write the convex differential inequality  at an arbitrary $\xi \in \mathcal H $ 
\begin{equation*}\label{Energ27}
\Phi(\xi) \geq \Phi (x(t)) +  \langle \nabla \Phi (x(t)), \xi - x(t)\rangle  .
\end{equation*}
Using that $x(\cdot)$ is bounded and (\ref{Energ25}) we immediately obtain
\begin{equation*}\label{Energ28}
\Phi(\xi) \geq \limsup_{t \to +\infty} \Phi (x(t)) \geq \liminf_{t \to +\infty} \Phi (x(t)) \geq \inf_{\mathcal H} \Phi .
\end{equation*}
This inequality being true for any $\xi \in \mathcal H $, we obtain (\ref{Energ26}).

%\medskip

\subsubsection*{Step 2: Using distance to equilibria as a Lyapunov function}

For an element $z\in S=\argmin\{ \Psi|\argmin\Phi\}$, we define the function $h_z:\R_+\to \R_+$
  by
$$
h_z(t)=\frac{1}{2}\|x(t) - z \| ^2.
$$
First write the optimality condition for $z \in S$. Recall that $C =  \argmin\Phi$.  Since $ z\in   \argmin (\Psi + \delta_C)$, we have 
$0 \in  \nabla \Psi (z) + \partial \delta_C (z)$. Since  $\delta_C (z) =  N_C (z)$, 
\begin{equation*}\label{Liapdist1}
 0 \in  \nabla \Psi (z) + N_C (z) .
\end{equation*}
Equivalently,
\begin{equation}\label{Liapdist2}
\exists p \in  N_C (z) \   \mbox{ such that } \   -p =  \nabla \Psi (z).
\end{equation}
Using $h_z$ as a Lyapunov-like function provides further informations that are described in the following lemma.
\begin{lemma}\label{Liapdist}  Let us assume that standing assumptions $(\mathcal H_0)$ hold together with the moderate decrease properties $(\mathcal H_1)$  and $(\mathcal H_2)$. 
Then, for any trajectory $x(\cdot):  [0, +\infty[ \to \mathcal H $  of \eqref{MIG}, the following properties hold:
\begin{itemize}
\item (i) for 	any $z \in S$,   \      $\lim_{t \to  +\infty} h_z(t)$ \ exists;
\item (ii) $\int_0^{+\infty} \Phi (x(t)) dt < + \infty$;
\item (iii) $\liminf_{t \to  +\infty} \Psi(x(t)) \geq \Psi (z)$.
\end{itemize} 
\end{lemma}
\textbf{Proof of Lemma \ref{Liapdist}}.
Given $z \in S$, let us compute the first and second order time derivatives of $h_z(\cdot)$. Using the classical derivation chain rule, we obtain
\begin{eqnarray}\label{deriv}
\dot{h}_z(t) &=& \langle  x(t)-z, \dot{x}(t) \rangle  \\
\ddot{h}_z(t) &=& \langle  x(t)-z, \ddot{x}(t) \rangle +   \|\dot{x}(t) \|^2.
\end{eqnarray}
By using \eqref{MIG} equation it follows that
\begin{eqnarray*}
\ddot{h}_z(t) + \gamma  \dot{h}_z(t) &=& \langle  x(t)-z, \ddot{x}(t) +   \gamma \dot{x}(t)  \rangle  +   \|\dot{x}(t) \|^2  \\
&=& \langle  x(t)-z, - \nabla \Phi(x(t)) -  \epsilon (t)  \nabla \Psi(x(t))\rangle +   \|\dot{x}(t) \|^2.
\end{eqnarray*}
Equivalently
\begin{equation}\label{Liapdist3}
\ddot{h}_z(t) + \gamma  \dot{h}_z(t) + \langle   \nabla \Phi(x(t)), x(t)-z \rangle   + \epsilon (t) \langle   \nabla \Psi(x(t)), x(t)-z   \rangle =  \|\dot{x}(t) \|^2  .
\end{equation}
Let us now write the convex differential inequalities (recall  that $z\in\argmin\Phi$, and hence $ \Phi (z) = 0$)
\begin{equation}\label{Liapdist4}
0 = \Phi (z) \geq  \Phi(x(t)) + \langle   \nabla \Phi(x(t)), z - x(t) \rangle,
\end{equation}
\begin{equation}\label{Liapdist5}
\Psi (z) \geq  \Psi(x(t)) + \langle   \nabla \Psi(x(t)), z - x(t) \rangle.
\end{equation}
Combining (\ref{Liapdist3}), (\ref{Liapdist4}) and (\ref{Liapdist5}), we obtain 
\begin{equation}\label{Liapdist6}
\ddot{h}_z(t) + \gamma  \dot{h}_z(t) + \Phi(x(t))   + \epsilon (t)\left( \Psi(x(t)) - \Psi(z) \right) \leq  \|\dot{x}(t) \|^2  .
\end{equation}
Let us now use the optimality condition (\ref{Liapdist2}), i.e., $  -p =  \nabla \Psi (z)$  with $p \in  N_C (z)$. We have 
\begin{equation}\label{Liapdist7}
\Psi(x(t))    \geq  \Psi (z)  + \langle   -p,  x(t) - z \rangle.
\end{equation}
Combining (\ref{Liapdist6}) and (\ref{Liapdist7})
\begin{equation}\label{Liapdist8}
\ddot{h}_z(t) + \gamma  \dot{h}_z(t)   \leq  \left(\epsilon (t) \langle  p,  x(t) \rangle  - \Phi(x(t))\right)    - \epsilon (t)  \langle  p,  z \rangle  +  \|\dot{x}(t) \|^2  .
\end{equation}
By definition of the Fenchel conjugate
\begin{equation}\label{Liapdist9}
\epsilon (t) \langle  p,  x(t) \rangle  - \Phi(x(t))  \leq  \Phi^* \left(\epsilon (t) p \right).
\end{equation}
On the other hand, by $p \in  N_C (z)$, we have
\begin{equation}\label{Liapdist10}
\epsilon (t)  \langle  p,  z \rangle =  \sigma_C\left(\epsilon (t) p\right). 
\end{equation}
Combining (\ref{Liapdist9}), (\ref{Liapdist10}) with  (\ref{Liapdist8})  we finally obtain
\begin{equation}\label{Liapdist11}
\ddot{h}_z(t) + \gamma  \dot{h}_z(t)   \leq   \left[ \Phi^* \left(\epsilon (t) p \right) - \sigma_C \left(\epsilon (t) p\right) \right]   +  \|\dot{x}(t) \|^2  .
\end{equation}
By  (\ref{Liapdist2}),   $p$ belongs to the cone $\mathcal C$ as defined in $(\mathcal H_2)$. 
Hence, by assumption $(\mathcal H_2)$ we have that 
\begin{equation*}\label{Liapdist12}  
\int_{0}^{+\infty}  \left[\Phi^* \left(\epsilon (t) p\right) - \sigma_C\left(\epsilon (t) p\right)\right]dt < + \infty.
\end{equation*}
On the other hand, we know by Lemma \ref{energy1} that $\int_{0}^{+\infty} \|\dot{x}(t) \|^2 dt < + \infty$.\\
Hence the second member of (\ref{Liapdist11}) is integrable on $[0, +\infty [$, and the  convergence of $h_z(t)$ as $t\to +\infty$ is a consequence of the following lemma, see \cite[Lemma 2.2 and 2.3]{Al}, \cite[Lemma 3.1]{AtCz1}.
\begin{lemma}\label{Alvlemma}  Let $t_0 \in \mathbb R$ and $h \in \mathcal C^2([t_0, + \infty [, \mathbb R^+)$ satisfy the following differential inequality
$$    \ddot{h}(t)  + \gamma \dot{h}(t)  \leq g(t)$$ 
with $g \in L^1([t_0, + \infty [, \mathbb R^+)$. Then $(\dot{h})_{+}$, the positive part of $\dot{h}$, belongs to $L^1([t_0, + \infty [, \mathbb R)$, and, as a consequence,
$\lim_{t \to  +\infty} h(t)$ exists.
\end{lemma}
\noindent Thus item $(i)$ of Lemma \ref{Liapdist} is proved. In order to obtain further estimates, let us return to (\ref{Liapdist6}). The same computation as above relying on the definition of the Fenchel conjugate yields
\begin{equation}\label{Liapdist13}  
 -\left[\Phi^* \left(\epsilon (t) p\right) - \sigma_C\left(\epsilon (t) p\right)\right] \leq \Phi(x(t))   + \epsilon (t)\left( \Psi(x(t)) - \Psi(z) \right).
\end{equation}
Combining (\ref{Liapdist6}) with (\ref{Liapdist13}), we obtain
\begin{equation*}\label{Liapdist14}  
\ddot{h}_z(t) + \gamma  \dot{h}_z(t) -\left[\Phi^* \left(\epsilon (t) p\right) - \sigma_C\left(\epsilon (t) p\right)\right] \leq \ddot{h}_z(t) + \gamma  \dot{h}_z(t) + \Phi(x(t))   + \epsilon (t)\left( \Psi(x(t)) - \Psi(z) \right) \leq  \|\dot{x}(t) \|^2.
\end{equation*}
 Let us  integrate the above inequalities between two large real numbers. Using that $\lim_{t \to  +\infty} h(t)$ exists, as it has been just proved,  noticing $\lim_{t \to  +\infty} \dot{h}_z(t) = 0$
(this last point results from (\ref{deriv}), the convergence of $\dot{x}(t)$ to zero, and the boundedness of $x(\cdot)$), and using assumption $(\mathcal H_2)$, we conclude that 
\begin{equation}\label{Liapdist15}  
\lim_{t \to  +\infty} \int_0^t \left[ \Phi(x(s))   + \epsilon (s)\left( \Psi(x(s)) - \Psi(z) \right)\right]ds   \quad  \mbox{exists}.
\end{equation}
Let us now use the same device as  \cite{AtCz2}, and split in equation (\ref{Liapdist8}) the term $\Phi(x(t))$ into two parts, so obtaining
\begin{eqnarray}\label{Liapdist106} 
\ddot{h}_z(t) + \gamma  \dot{h}_z(t) + \frac{1}{2} \Phi(x(t))  &\leq &   \|\dot{x}(t) \|^2   +  \frac{1}{2} \left[  \langle 2 \epsilon (t) p,  x(t) \rangle  - \Phi(x(t))\right]   
 -\frac{1}{2}  \langle 2\epsilon (t)  p,  z \rangle  \\
& \leq &  \|\dot{x}(t) \|^2   +  \frac{1}{2}\left[ \Phi^* \left(2\epsilon (t) p\right) - \sigma_C\left(2\epsilon (t) p\right)\right]. \nonumber
\end{eqnarray}
Since $\mathcal C$ is a cone, $2p$ still belongs to $\mathcal C$, which by assumption $(\mathcal H_2)$ yields integrability of the second member of inequality (\ref{Liapdist106}).
It follows at once that 
\begin{equation}\label{Liapdist17}
\int_0^{+\infty} \Phi (x(t)) dt < + \infty,
\end{equation}
which proves item $(ii)$.
Combining (\ref{Liapdist15})  and (\ref{Liapdist17}), we also get
\begin{equation*}\label{Liapdist18}
\lim_{t \to  +\infty} \int_0^t  \epsilon (s)\left( \Psi(x(s)) - \Psi(z) \right)ds   \quad  \mbox{exists}.
\end{equation*}
Using assumption $(\mathcal H_1)$, i.e., $\int_0^{+\infty} \epsilon (t) dt = + \infty$, we obtain item $(iii)$, that is
\begin{equation*}\label{Liapdist19}
\liminf_{t \to  +\infty}  \Psi(x(t)) \leq  \Psi(z).
\end{equation*}

%\medskip

\subsubsection*{Step 3: Using a rescaled energy as a Lyapunov function}

\noindent Detecting the property of asymptotic hierarchical minimization requires  further Lyapunov analysis based on the rescaled energy function

\begin{equation*}\label{Energ2}
E_2(t) =   \frac{1}{2\epsilon(t)} \|\dot{x}(t)\|^2   +   \frac{1}{\epsilon(t)}  \Phi(x(t))   + \Psi (x(t)).
\end{equation*}
\begin{lemma}\label{Liapresc}  Let us assume that standing assumptions $(\mathcal H_0)$ hold together with the moderate decrease properties $(\mathcal H_1)$, $(\mathcal H_2)$ and $(\mathcal H_3)$. 
Then, for any trajectory $x(\cdot):  [0, +\infty [ \to \mathcal H $  of \eqref{MIG}, and any $z\in S$,
the following properties hold:
\begin{itemize}
\item (i) \  $  \dot{E_2} (t) =   - \frac{\gamma}{\epsilon(t)} \|\dot{x}(t)\|^2  -  \frac{\dot{\epsilon}(t)}{{\epsilon(t)}^2} \left( \frac{1}{2} \|\dot{x}(t)\|^2   +   \Phi(x(t))   \right)$;            \item (ii) \  $\lim_{t \to  +\infty} E_2(t)$  \  exists           ;
\item (iii) \   $\lim_{t \to  +\infty} \Psi(x(t)) = \Psi (z)$.
\end{itemize} 
\end{lemma}
\textbf{Proof of Lemma \ref{Liapresc}}.
Let us compute the first  order time derivative of $E_2(\cdot)$. By using classical derivation chain rule, and equation \eqref{MIG}, we obtain
\begin{eqnarray*}
\dot{E_2}(t)   & = &  \frac{1}{\epsilon(t)} \langle \dot{x}(t), \ddot{x}(t)\rangle -   \frac{\dot{\epsilon}(t)}{2{\epsilon(t)}^2}\|\dot{x}(t)\|^2 
+ \frac{1}{\epsilon(t)} \langle \nabla \Phi (x(t)), \dot{x}(t)\rangle  -   \frac{\dot{\epsilon}(t)}{{\epsilon(t)}^2} \Phi(x(t))  + \langle \nabla \Psi (x(t)), \dot{x}(t)\rangle   \\
& =  &  \frac{1}{\epsilon(t)}  \langle \dot{x}(t), \ddot{x}(t) +  \nabla \Phi (x(t)) + \epsilon(t) \nabla \Psi (x(t)) \rangle  
 -   \frac{\dot{\epsilon}(t)}{{\epsilon(t)}^2} \left(\frac{1}{2}  \|\dot{x}(t)\|^2      + \Phi (x(t))  \right)\\
& =  &      - \frac{\gamma}{\epsilon(t)} \|\dot{x}(t)\|^2  -  \frac{\dot{\epsilon}(t)}{{\epsilon(t)}^2} \left( \frac{1}{2} \|\dot{x}(t)\|^2   +   \Phi(x(t))   \right),      
\end{eqnarray*}
which proves item $(i)$.
By assumption $(\mathcal H_3)$, there exists some positive constant $k$ such that 
\begin{equation*}
0 \leq - \dot{\epsilon}(t) \leq k \epsilon^2 (t).  
\end{equation*}
As a consequence
\begin{equation*}
\dot{E_2}(t)   \leq  k  \left( \frac{1}{2} \|\dot{x}(t)\|^2   +   \Phi(x(t))   \right).
\end{equation*}
Relying on the integrability property on $]0, + \infty[$   of $ \|\dot{x}(\cdot)\|^2 $  (see Lemma \ref{energy1} $(ii)$),  and  $ \Phi(x(\cdot))$  (see Lemma \ref{Liapdist}  $(ii)$),
we deduce that $\dot{E_2}(\cdot)$ is majorized by an integrable function. Since $E_2$ is bounded from below, it follows that
\begin{equation*}
\lim_{t \to  +\infty} E_2(t)  \  \mbox{exists}, 
\end{equation*}
that's item $(ii)$. Since $E_2(t) \geq \Psi (x(t))$, by using Lemma \ref{Liapdist} $(iii)$ we obtain
\begin{eqnarray*}
\lim_{t \to  +\infty} E_2(t)  &\geq & \liminf_{t \to  +\infty} \Psi(x(t))\\
                             & \geq &  \Psi (z).
\end{eqnarray*}
Let us prove that 
\begin{equation*}
\lim_{t \to  +\infty} E_2(t)  = \Psi (z).
\end{equation*}
To achieve this, we argue by contradiction and assume that
\begin{equation*}
\lim_{t \to  +\infty} E_2(t)  > \Psi (z).
\end{equation*}
It will follow that for $t$ sufficiently large  and for some positive $\alpha$
\begin{equation*}
\frac{1}{2\epsilon(t)} \|\dot{x}(t)\|^2   +   \frac{1}{\epsilon(t)}  \Phi(x(t))   + \Psi (x(t)) > \Psi (z) + \alpha.
\end{equation*}
Hence
\begin{equation*}
 \frac{1}{2} \|\dot{x}(t)\|^2   +   \Phi(x(t))   +   \epsilon(t)   \left( \Psi (x(t)) - \Psi (z) \right) >  \alpha \epsilon(t).
\end{equation*}
By integration of this inequality, using the fact that  $\int_0^{+\infty} \epsilon (t) dt = + \infty$ (assumption $(\mathcal H_1)$), and the fact that (equation (\ref{Liapdist15}))
\begin{equation*}
\lim_{t \to  +\infty} \int_0^t \left[ \Phi(x(s))   + \epsilon (s)\left( \Psi(x(s)) - \Psi(z) \right)\right]ds   \quad  \mbox{exists}
\end{equation*}
we obtain a contradiction.
Thus 
 \begin{equation}\label{energy3}
\Psi (z)  =  \lim_{t \to  +\infty} E_2(t)  \geq  \limsup_{t \to  +\infty} \Psi (x(t)).
\end{equation}
Combining (\ref{energy3}) with  Lemma \ref{Liapdist} $(iii)$, we finally obtain item $(iii)$
\begin{equation*}\label{energy4}
\lim_{t \to  +\infty} \Psi(x(t)) = \Psi (z).
\end{equation*}

\noindent \textit{End of the proof of Theorem ~\ref{basicthm}}. 

\noindent We can now apply Opial's lemma.  By Lemma \ref{Liapdist} $(i)$,  for
every $z\in S$ there is convergence of $ \|x(.) -z\|$. By Lemma  \ref{energy1} $(iv)$
$\lim_{t \to  +\infty}  \Phi(x(t))= \inf_{\mathcal H} \Phi = 0$. This clearly implies that every weak cluster point of $x(\cdot)$ belongs to
$C$. Moreover  by Lemma \ref{Liapresc} $(iii)$, $\lim_{t \to  +\infty} \Psi(x(t)) = \Psi (z)$. This forces every weak cluster point of $x(\cdot)$ to belong to $S$.
In the  two preceding results, we use the  lower semicontinuity of $\Phi$ and $\Psi$ for the weak topology. $\Box$

\section{Strong convergence results}\label{section_strong}

The  \eqref{MIG} system was first introduced in the strongly convex case by~\cite{AtCz1} and its generalization~\cite{Cabot}. The proof relies on a geometrical argument which uses the weak compactness of the set:
$$\{x | \langle \nabla \Psi(x), x-x_\infty\rangle\leq 0\},$$
where $x_\infty$ is the unique point in $S$. % A straightforward computation yields the weak compacity of the set $\{x | \langle \nabla \Psi(x), x-z\rangle\leq 0\}$ under the assumption of strong convexity of the potential $\Psi$ in it classical form:
%$$\exits K>0, \forall x\in H,\forall y\in H, \langle \nabla \Psi(x)- \nabla \Psi(y), x-y\rangle\geq K   |x-y|^2.$$
%However, it is not clear that the set is weakly compact under the assumption of strong convexity of the potential $\Psi$ in it weak form in~\cite{Cabot}: 
%$$\forall x\in H,\forall y\in H, \langle \nabla \Psi(x)- \nabla \Psi(y), x-y\rangle\geq f( |x-y|),$$
%for some function $f:\R_+\to\R_+$ such that $f(t_n)\to 0 \Rightarrow t_n \to 0$.
%In this case, o
Our approach provides a different hindsight of the convergence result, as follows.

 \begin{Theorem} \label{thm_strong}
Let us suppose the standing assumptions $(\mathcal H_0)$, and $ (\mathcal H_1)  \   \mbox{ }\mbox{ } \mbox{ }    \int_0^{+\infty} \epsilon (t) dt = + \infty$. Assume that the potential  $\Psi$ is uniformly convex, i.e., there exists a function $\theta:\R_+\to\R_+$ such that $\theta(t_n)\to 0 \Rightarrow t_n \to 0$ and, for every $x$ and $y$ in $H$, 
$$\langle \nabla \Psi(x)- \nabla \Psi(y), x-y\rangle\geq \theta({ \|x-y\|}).$$
Let $x_\infty$ be the unique point in  $S=\argmin\{
  \Psi|\argmin\Phi\}$. Assume one of the following:\\
$(i)$ \cite[Theorem 2.3]{AtCz1}, \cite[Theorem 5.1]{Cabot} the set $\{x | \langle \nabla \Psi(x), x- {x_\infty}\rangle\leq 0\}$ is weakly compact.\\
$(ii)$ \begin{align*}
	&(\mathcal H_2) \   \mbox{ }\mbox{ } \mbox{ }  \forall p\in \mathcal C \mbox{ } \mbox{ }\mbox{ } \int_{0}^{+\infty}  \left[\Phi^* \left(\epsilon (t) p\right) - \sigma_C\left(\epsilon (t) p\right)\right]dt 
< + \infty.\\
 & (\mathcal H_3)  \     \mbox{ }\mbox{ } \mbox{ }    0 \leq - \dot{\epsilon}(t) \leq k \epsilon^2 (t)  \  \mbox{ }\mbox{ }  \mbox{for some positive constant} \  k, \mbox{ }  \mbox{and large} \mbox{ } t.
\end{align*}
Let $x$ be a classical maximal solution of $\eqref{MIG}$. Then $\lim_{t\to +\infty}{ \|x(t)-x_\infty \|}=0$
\end{Theorem}

%\begin{remark} As stated above, if for some  $ K>0$, and every $x$ and $y$ in $H$,
%$$ \langle \nabla \Psi(x)- \nabla \Psi(y), x-y\rangle\geq K   |x-y|^2, $$
%the set 
%$\{x | \langle \nabla \Psi(x), x-z\rangle\leq 0\}$ is weakly compact. Indeed, write
%$$\langle \nabla \Psi(x), x-z\rangle= \langle \nabla \Psi(x)- \nabla \Psi(z), x-z\rangle+\langle \nabla \nabla \Psi(z), x-z\rangle, $$
%thus
%$$\langle \nabla \Psi(x), x-z\rangle\geq K|x-z|^2-\frac{K}{2}|x-z|^2-\frac{2K}|\nabla \Psi(z)$$

 \textbf{Proof of Theorem~\ref{thm_strong}.}\\
\textbf{Case (i)}
We refer to~\cite{AtCz1} and its generalization~\cite{Cabot}.\\
\textbf{Case (ii)} Just follow the proof of Theorem~\ref{basicthm}. With $x_\infty$ the unique point in $S$, 
Equations (\ref{Liapdist3}) and (\ref{Liapdist4}), together with the strong convexity of $\Psi$, yield
$$
\ddot{h}_{x_\infty  }(t) + \gamma  \dot{h}_{x_\infty  }(t) + \Phi(x(t))   + \epsilon (t)\langle   \nabla \Psi(x(t)), x(t)-x_\infty   \rangle \leq  \|\dot{x}(t) \|^2  .
$$
Writing
$$\langle \nabla \Psi(x), x-{x_\infty}\rangle= \langle \nabla \Psi(x)- \nabla \Psi({x_\infty}), x-{x_\infty}\rangle+\langle {\nabla} \Psi({x_\infty}), x-{x_\infty}\rangle, $$
and following the proof of  Lemma \ref{Liapdist}, we deduce
$$
\ddot{h}_{x_\infty  }(t) + \gamma  \dot{h}_{x_\infty  }(t) +\epsilon(t)\theta({\|x(t)-x_\infty  \|})  \leq   \left[ \Phi^* \left(\epsilon (t) p \right) - \sigma_C \left(\epsilon (t) p\right) \right]   +  \|\dot{x}(t) \|^2  ,
$$
with
$  p = - \nabla \Psi (x_\infty )$  with $p \in  N_C (x_\infty )$. Thus
$$
\int_{0}^{+\infty} \epsilon(t) \theta({\|x(t)-x_\infty  \|})<+\infty.
$$
Once the convergence of ${h}_{x_\infty  }$ is obtained, its limit can only be $0$.

\section{Multiscale aspects}\label{multiscale}

\subsection{Invariance by affine rescaling} Let us show that the moderate decrease conditions are invariant by affine rescaling. In particular this means that our asymptotic analysis of \eqref{MIG} system 
is not affected by an affine change of variable. Given $t \mapsto x(t)$  a solution trajectory of \eqref{MIG}, and $a >0$ a positive parameter (coefficient of the affine rescaling) let us set
 \begin{equation*}\label{rs1}
y(t) := x(at).
\end{equation*}
Direct application of the  derivation chain rule yields
 \begin{equation*}\label{rs2}
 \qquad     \ddot{y}(t)   +  a \gamma  \dot{y}(t) +    a^2 \nabla \Phi (y(t)) +   a^2 \epsilon (at) \nabla \Psi  (y(t)) = 0.
\end{equation*}
This is still a \eqref{MIG} system with the new data: 
\begin{align*}
\tilde{\Phi}(y) & =  a^2 \Phi (y)\\
\tilde{\Psi}(y) & =  \Psi (y)\\
\tilde{\epsilon}(t) & =  a^2 \epsilon (at)
\end{align*}
Since $\argmin \tilde{\Phi} = \argmin \Phi = C$, the cone $\mathcal C$ is invariant by affine rescaling, and  condition $(\mathcal H_2)$ for the rescaled equation (\ref{rs2}) is equivalent to 
 \begin{equation*}\label{rs3}
   \mbox{ }\mbox{ } \mbox{ }  \forall p\in \mathcal C \mbox{ } \mbox{ }\mbox{ } \int_{0}^{+\infty}  \left[\tilde{\Phi}^* \left(\tilde{\epsilon}(t) p\right) - \sigma_C\left(\tilde{\epsilon} (t) p\right)\right]dt 
< + \infty.
\end{equation*}
Since 
\begin{align*}
\tilde{\Phi}^* \left(\tilde{\epsilon}(t) p\right) & =  a^2 {\Phi}^* \left(\epsilon(at) p\right)  \\
\sigma_C\left(\tilde{\epsilon} (t)p\right) & =  a^2 \sigma_C\left(\epsilon(t)p\right)
\end{align*}
we have
\begin{align*}
 \int_{0}^{+\infty}  \left[\tilde{\Phi}^* \left(\tilde{\epsilon}(t) p\right) - \sigma_C\left(\tilde{\epsilon} (t) p\right)\right]dt & = a^2 \int_{0}^{+\infty}  \left[{\Phi}^* \left(\epsilon(at) p\right) - \sigma_C\left(\epsilon (at) p\right)\right]dt \\
 & = a \int_{0}^{+\infty}  \left[{\Phi}^* \left(\epsilon(t) p\right) - \sigma_C\left(\epsilon (t) p\right)\right]dt
< + \infty.
\end{align*}
Thus, we recover the $(\mathcal H_2)$ condition on the initial sytem. Clearly conditions $(\mathcal H_1)$ and $(\mathcal H_3)$ are also invariant by affine rescaling.

\subsection{Asymptotic penalization} 
In this section, by changing the time scale, we show dynamic systems that are connected and having multiscale aspects.
We start from 
\begin{equation}
\tag{\ref{MIG}} 
   m\ddot{x}(t)   + \gamma  \dot{x}(t) +    \nabla \Phi (x(t)) +   \epsilon (t) \nabla \Psi  (x(t)) = 0,
\end{equation}
which involves a  control term which is asymptotically vanishing ($\epsilon (t) \to 0$  as $t\to +\infty$). The nonnegative (mass) parameter $m$ has been introduced in the equation to cover both the second order and the first order system (obtained by taking $m=0$).
We use the following dictionary, introduced in \cite{AtCz2}.

\begin{lemma}[dictionary]\label{6_2_2009} Let $T_\beta$ and $T_\varepsilon$ be two
  elements in $\left(\R_+\setminus\{0\}\right)\cup{+\infty}$. Take two functions of
  class $C^1$
\begin{eqnarray*}
\beta : [0,T_\beta [&\to &\R_+\setminus\{0\};\\
\varepsilon: [0,T_\varepsilon [&\to &\R_+\setminus\{0\}.
\end{eqnarray*}
Define $ t_\beta:  [0,T_\varepsilon [\to [0,T_\beta [$  \ and \ $
    t_\varepsilon: [0,T_\beta [ \to [0,T_\varepsilon)$ by
$$
\int_0^{t_\beta(t)} \beta(s)ds =t \mbox{ and }
\int_0^{t_\varepsilon(t)} \varepsilon(s)ds =t
.
$$
Assume that, for every $t$, 
$$
\varepsilon(t)\beta(t_\beta(t))=1.
$$
Then 
\begin{eqnarray*}
t_\varepsilon \circ t_\beta={\rm id}_{[0,T_\varepsilon [} & and &
  T_\varepsilon=\int_0^{T_\beta} \beta;\\
t_\beta \circ t_\varepsilon={\rm id}_{[0,T_\beta [} & and &
  T_ \beta=\int_0^{T_\varepsilon}\varepsilon.
\end{eqnarray*}
If $x$ is a solution trajectory of 
\begin{equation}
\tag{\ref{MIG}} \qquad     m\ddot{x}(t)   + \gamma  \dot{x}(t) +    \nabla \Phi (x(t)) +   \epsilon (t) \nabla \Psi  (x(t)) = 0,
\end{equation}
then $w:=x\circ t_\varepsilon$ is a  solution trajectory of 
\begin{equation}\label{multi0}
\tag{{\rm MIG$_{\beta}$}}  \frac{m}{\beta (t)}\ddot{w}(t) + (\gamma  - m\frac{\dot{\beta}(t)}{{\beta}^2 (t)}  )\dot{w}(t) +  \beta(t) \nabla \Phi (w(t)) +            
\nabla \Psi (w(t))= 0.
\end{equation}
Conversely, if $w$  is a  solution of \eqref{multi0},
then $w\circ t_\beta$ is  a  solution of \eqref{MIG}.

\end{lemma}

{\bf Proof of Lemma~\ref{6_2_2009}.}   Let us make the change of variable associated with function $t_\varepsilon(\cdot)$. Note that $\dot{t}_\varepsilon(t)= \beta (t)$.
Let $x(\cdot)$ be  a  solution
of  \eqref{MIG} and write the system  \eqref{MIG} at the point ${t_\varepsilon(t)}$:
$$
 m\ddot{x}(t_\varepsilon(t))   +  \gamma \dot{x}(t_\varepsilon(t)) + 
  \nabla \Phi (x(t_\varepsilon(t))) +  \varepsilon(t_\varepsilon(t) ) \nabla
\Psi (x(t_\varepsilon(t)))= 0.
$$
After multiplication  by $\dot{t}_\varepsilon(t)$ we get
$$
 m\dot{t}_\varepsilon(t)\ddot{x}(t_\varepsilon(t))   +  \gamma \dot{t}_\varepsilon(t)\dot{x}(t_\varepsilon(t)) + 
 \dot{t}_\varepsilon(t) \nabla \Phi (x(t_\varepsilon(t))) + \dot{t}_\varepsilon(t) \varepsilon(t_\varepsilon(t) ) \nabla
\Psi (x(t_\varepsilon(t)))= 0.
$$
Set $w=x\circ t_\varepsilon$.
% The map $w=x\circ t_\varepsilon$ is absolutely continuous, 
According to $\dot{w}(t)= \dot{t}_{\varepsilon}(t)\dot{x}(t_{\varepsilon}(t))$,
$\dot{t}_\varepsilon(t)= \beta (t)$ and $ \dot{t}_\varepsilon(t) \varepsilon(t_\varepsilon(t) )=1$, we obtain
\begin{equation}\label{multi1}
 m\dot{t}_\varepsilon(t)\ddot{x}(t_\varepsilon(t))   + \gamma \dot{w}(t) +  \beta(t) \nabla \Phi (w(t)) +            
\nabla \Psi (w(t))= 0.
\end{equation}
From
$$
\ddot{w}(t) =  {\dot{t}_\varepsilon}^2(t)  \ddot{x}(t_\varepsilon(t)) + 
\ddot{t}_{\varepsilon} (t)\dot{x}(t_\varepsilon(t)),
$$
we obtain
\begin{equation}\label{multi2}
\dot{t}_\varepsilon(t)\ddot{x}(t_\varepsilon(t)) = \frac{1}{\beta (t)}\ddot{w}(t)  - \frac{\dot{\beta}(t)}{{\beta}^2 (t)}  )\dot{w}(t). 
\end{equation}
Combining (\ref{multi1}) and (\ref{multi2}) gives the result.
$\Box$

\medskip

Accordingly, all our results can be written for the system \eqref{multi0}.  Let us formulate the conditions of moderate decrease of  $\epsilon(\cdot)$ with the help of $\beta (\cdot)$.

\begin{align*}
	& (\mathcal H_1)  \   \mbox{ }\mbox{ } \mbox{ }    \int_0^{+\infty} \epsilon (t) dt = + \infty  \  \mbox{becomes}  \  \beta (t) \to + \infty \ \mbox{as}  \ t \to + \infty;
	\\
	&(\mathcal H_2) \   \mbox{ }\mbox{ } \mbox{ }  \forall p\in \mathcal C \mbox{ } \mbox{ }\mbox{ } \int_{0}^{+\infty}  \left[\Phi^* \left(\epsilon (t) p\right) - \sigma_C\left(\epsilon (t) p\right)\right]dt 
< + \infty \ \mbox{becomes} \  \int_{0}^{+\infty} \beta (t) \left[\Phi^* \left(\frac{p}{ \beta (t)}\right) -
\sigma_C \left(\frac{p}{ \beta (t)}\right)\right]dt < + \infty;  \\
 & (\mathcal H_3)  \     \mbox{ }\mbox{ } \mbox{ }    0 \leq - \dot{\epsilon}(t) \leq k \epsilon^2 (t)  \  \mbox{ }\mbox{ }  \mbox{for some positive constant} \  k, \mbox{ }  \mbox{and large} \mbox{ } t   \ \mbox{becomes} \   \dot{\beta}(t) \leq k \beta (t).
\end{align*}
Note that $(\mathcal H_1)$ and $(\mathcal H_3)$ imply $\lim_{t \to + \infty} \frac{\dot{\beta}(t)}{{\beta}^2 (t)}=0$. 
Thus \eqref{multi0} asymptotically behaves like
\begin{equation*}\label{multi3}
\qquad \frac{m}{\beta (t)}\ddot{w}(t) + \gamma  \dot{w}(t) +  \beta(t) \nabla \Phi (w(t)) +            \nabla \Psi (w(t))= 0.
\end{equation*}
Note that the coefficient $\frac{m}{\beta (t)}$ in front of the acceleration term goes to zero as time $t$ goes to infinity.
It is an interesting subject for further research to study the related system 
\begin{equation*}\label{multi4}
 \qquad \ddot{w}(t) + \gamma  \dot{w}(t) +  \beta(t) \nabla \Phi (w(t)) +            \nabla \Psi (w(t))= 0,
\end{equation*}
where the coefficient in front of the acceleration term is a fixed real positive number.

\if

Before doing so, let us analyse the
corresponding assumption $(\mathcal H_1)$.

\begin{remark}\label{r1}{About Assumption ($\mathcal H_1$).} Take $\beta(\cdot)$ and
  $\varepsilon(\cdot)$ as in Lemma~\ref{6_2_2009}, with
  $T_\beta=T_\varepsilon=+\infty$.
  Then, by making the change of variable $t =  t_\beta(s)$ in the following integral, and by using $\dot{t}_\beta(t)\beta(t_\beta(t))=1$, we obtain

\begin{eqnarray*}
\int_{0}^{+\infty} \beta (t) \left(\Psi^*  \left(\frac{p}{ \beta (t)}\right) -
\sigma_C  \left(\frac{p}{ \beta (t)}\right)\right)dt & = &  \int_{0}^{+\infty} \dot{t}_\beta(s)\beta (t_\beta(s)) \left(\Psi^* \left(\frac{p}{ \beta (t_\beta(s))}\right) -
\sigma_C \left(\frac{p}{ \beta (t_\beta(s))}\right)\right)ds\\
&=&\int_{0}^{+\infty} \Psi^*
(\varepsilon(s)p ) - \sigma_C (\varepsilon(s)p)ds.
\end{eqnarray*}

\noindent Thus condition ($\mathcal H_1$) becomes

$$\forall p\in R(N_C) \ \ \int_{0}^{+\infty} \Psi^*
(\varepsilon(t)p ) - \sigma_C (\varepsilon(t)p)dt < + \infty.$$

\noindent In  particular, when  $\Psi (x) =  \frac{1}{2} dist^2 (x, C)$,
then $\Psi^* (x) =  \frac{1}{2} \|x\|^2  +  \sigma_C (x)$ and 
 
$$(\mathcal H_1) \Longleftrightarrow  \int_{0}^{+\infty} \frac
{1}{\beta (t)} dt < + \infty\Longleftrightarrow
\int_{0}^{+\infty}\varepsilon(t)^2 dt < + \infty. \ \Box$$ 
\end{remark}
 
\begin{Theorem} \label{Theweak_eps}
Let  $\Phi$ and $\Psi$ satisfy the assumptions of
Theorem ???. Let us assume that,

\begin{itemize}
\item
${(\mathcal H_{1})}_ {\varepsilon}$\mbox{ }\mbox{ }  $ \forall p\in R(N_C),$
  $\displaystyle \int_{0}^{+\infty} \Psi^*
(\varepsilon(t)p ) - \sigma_C (\varepsilon(t)p)dt < + \infty.$

\item
${(\mathcal H_{2})}_ {\varepsilon}$\mbox{ }\mbox{ } $\varepsilon(\cdot)$ is
 a non increasing function of class $C^1$,
  such that $\lim_{t\to+\infty} \varepsilon (t)=0$, $\displaystyle \int_{0}^{+\infty}  \varepsilon (t)dt=+\infty$,
and for some $k\geq 0$, $-k \varepsilon^2\leq 
  \dot{\varepsilon}$.

\end{itemize}

Let $x(\cdot)$ be a strong solution of $(MAG)_{\varepsilon}$. Then:

\begin{eqnarray*}
(i)&\mbox{ weak convergence }&\exists x_{\infty}\in S=\argmin\{
  \Phi|\argmin\Psi\},\qquad  w-\lim_{t\to+\infty}
  x(t)=x_{\infty};\\
(ii)&\mbox{ minimizing properties } &\lim_{t\to+\infty} \Psi(x(t))=
  0;\\
&&\lim_{t\to +\infty} \Phi(x(t)) = \min\Phi |_{\argmin\Psi};\\
(iii)&&\forall z\in S  \lim_{t\to+\infty} \|x(t) -z\| \mbox{  exists  };\\
(iv)&\mbox{ estimations }&\lim_{t\to +\infty} \frac{1}{\varepsilon(t)}\Psi(x(t))=0;\\
&&  \int_{0}^{+\infty} \Psi (x(t)) dt < +\infty;\\
&&\limsup_{\tau\to +\infty} \int_{0}^{\tau}\varepsilon(t)\left(\Phi (x(t))- \min\Phi |_{\argmin\Psi}\right)  dt < +\infty.
\end{eqnarray*}
\end{Theorem}

\textbf{{Proof of Theorem} \ref{Theweak_eps}}
Equivalence between Theorem~\ref{Theweak_eps}
and  Theorem ??? is a 
consequence of Lemma~\ref{6_2_2009}, Remark~\ref{r1} and the equivalent formulation of condition $(\mathcal H_{2})$ as given below:

Write $(\mathcal H_{2})$ at the point $t_{\beta}(t)$ 

$$
\dot{\beta}(t_{\beta}(t))\leq k \beta(t_{\beta}(t))
$$

and multiply by $\dot{t_{\beta}}(t)$ (which is nonnegative)

$$
\dot{\beta}(t_{\beta}(t))\dot{t_{\beta}}(t) \leq k \beta(t_{\beta}(t))\dot{t_{\beta}}(t).
$$
  
Owing to $\beta(t_{\beta}(t))\dot{t_{\beta}}(t) = 1$, we get

$$
\frac{d}{dt} \beta(t_{\beta}(t)) \leq k,
$$

which, by using $\varepsilon(t)\beta(t_\beta(t))=1$, finally yields

$$
\frac{d}{dt} \left(\frac{1}{\varepsilon(t)} \right) =  - \frac{\dot{\varepsilon}(t)}{{\varepsilon}^2(t)}  \leq k
$$

that's  ${(\mathcal H_{2})}_ {\varepsilon}$. $\Box$

\fi

\section{Coupled gradient dynamics}\label{coupled_section}

Throughout this section we make the following assumptions: 
\begin{itemize}
\item $\mathcal H = \Xiu \times \Xid$ is the cartesian product of two Hilbert spaces, set $x=(x_1,x_2)$;
\item  $\Phi(x) = f_1(x_1) + f_2(x_2) $,  \  $f_1\in \Gamma_0 (\Xiu)$,  $f_2\in \Gamma_0 (\Xid)$ 
are convex functions which are continuously differentiable,  $\nabla f_i$ is Lipschitz continuous on bounded sets of $\X_i$, $i=1,2$.
We assume that $C_i = \argmin f_i$ is nonempty.
\item $\Psi (x) = \frac{1}{2}\|L_1x_1 - L_2x_2\|^2_{\Z}$,  \   
$L_1 \in L(\Xiu, \Z)$ and  $L_2 \in L(\Xid, \Z)$ are linear continuous operators acting respectively from $\Xiu$  and $\Xid$ into a third Hilbert space $\Z$;
\item $\epsilon : \R^+ \rightarrow \R^+$ is  a function of $t$ which tends to $0$ as $t$ goes to $+\infty$.
\end{itemize}
\noindent In this setting
\begin{equation}
\tag{\ref{MIG}}
    \ddot{x}(t)   + \gamma  \dot{x}(t) +    \nabla \Phi (x(t)) +   \epsilon (t) \nabla \Psi  (x(t)) = 0, 
\end{equation}
\noindent becomes
\begin{equation}\label{couplsys}
\left\{
\begin{array}{l}
\ddot{x_1}(t) + \gamma \dot{x_1}(t) +  \nabla f_1(x_1(t)) +  \epsilon (t)L_1^{*} ( L_1x_1(t) - L_2x_2(t))= 0\\
\rule{0pt}{25pt}
\ddot{x_2}(t) +  \gamma \dot{x_2}(t) + \nabla f_2(x_2(t))  + \epsilon (t) L_2^{*} ( L_2x_2(t) - L_1x_1(t))= 0.
\end{array}\right.
\end{equation}
\noindent Because of the quadratic property of $\Psi$, condition $({\mathcal H}_1)$ can be equivalently written 
\begin{equation}\label{couplsys1}
 \epsilon (\cdot) \in L^2 (0, + \infty) \setminus L^1 (0, + \infty).
\end{equation}
\noindent As a straight application of Theorem \ref{basicthm}, assuming (\ref{couplsys1}) and the growth condition 
\begin{equation*}\label{couplsys2}
  0 \leq - \dot{\epsilon}(t) \leq k \epsilon^2 (t)  \  \mbox{ }\mbox{ }  \mbox{for some positive constant} \  k, \mbox{ }  \mbox{and large} \mbox{ } t,
\end{equation*}
\noindent we obtain that
$x(t)= (x_1(t),x_2(t))\rightarrow x_{\infty}=(x_{1,\infty},x_{2,\infty})$ weakly in  $\mathcal H$
where $(x_{1,\infty},x_{2,\infty})$ is a solution of 
\begin{equation*}\label{structopt}
\mbox{min} \left\{ \|L_1x_1 - L_2x_2\|: \mbox{ }\mbox{ } x_1 \in \argmin f_i, \mbox{ } x_2 \in \argmin f_2\right\},
\end{equation*}
assuming that the solution set of the above problem is nonempty.
The above model describes two oscillators coupled by a force $\epsilon (t) \nabla \Psi  (x(t))$ which asymptotically vanishes. When the  control variable  $\epsilon (t)$  tends to zero in  a \emph{moderate} way, i.e., $\epsilon (\cdot) \in L^2 (0, + \infty) \setminus L^1 (0, + \infty)$, the control forces the two oscillators to stabilize asymptotically in equilibria  that are "as close as possible" to each other.

\smallskip

This model has a natural interpretation as a  Nash equilibration process for 
potential games. Consider  two  players with respective individual 
payoff functions $f_1 \in \Gamma_0 (\Xiu)$ and $f_2 \in  \Gamma_0 (\Xid)$.
Suppose that, at time $t$, the two players have a joint payoff 
$\frac{\epsilon (t)}{2} \| L_1x_1 - L_2x_2\|^2$. This attractive potential
may reflect some joint constraint on resources.
System (\ref{couplsys}) describes a continuous dynamic which, at time $t$, is attached to the potential team game (Monderer-Shapley, 1996) associated to the payoff functions
$$\left\{
\begin{array}{l}
F_1(x_1 ,x_2) = f_1 (x_1) +   \frac{\epsilon (t)}{2} \| L_1x_1 - L_2x_2\|^2                         \\ 
\rule{0pt}{20pt}
F_2(x_1 ,x_2) = f_2 (x_2) + \frac{\epsilon (t)}{2} \| L_1x_1 - L_2x_2\|^2    
\end{array}\right.$$
The second order nature of the dynamics, and the damping term reflects the importance of the inertia, and friction aspects in the process.
This model provides an explanation of what seems strange phenomenum, why specific equilibria are observed in real life, despites the fact that no force is present who can explain. It is history, and the fact that these forces disappear asymptotically, which gives a response. 

\section{Coupled wave equations}\label{coupledwavesection}
Consider the following infinite dimensional version of the situation discussed above.  It modelizes two  wave equations,  with Neumann boundary conditions, which are coupled by an asymptotic vanishing term. Let us give

\begin{itemize}
\item $\Omega$  a bounded open set in $\mathbb R^n$;

\smallskip

\item $\gamma >0$ (damping parameter), $\alpha_1 >0, \alpha_2 >0$  (waves propagation
speed) which are  positive constants;

\smallskip

\item $h_1, h_2 \in 
L^2 (\Omega)$ verifying $\int_{\Omega}h_1= \int_{\Omega}h_2 = 0$ (compatibility condition);

\smallskip

\item $u_{01}, u_{02} \in H^1 (\Omega)$, $v_{01}, v_{02} \in L^2 (\Omega)$ (Cauchy data);

\smallskip

\item $\epsilon : [0, +\infty[   \to {\mathbb R}_+$  a non-increasing function of class $\mathcal C^2$ such that $\lim_{t\to +\infty} \epsilon (t) =        \lim_{t\to +\infty} \dot{\epsilon} (t)=0.$ 
\end{itemize}

\smallskip

The  system  is written as follows: 

\begin{equation}\label{couplwave1}
\left\{
\begin{array}{l}
\frac{\partial^2 u_1}{\partial t^2} + \gamma \frac{\partial u_1}{\partial t}  - \alpha_1 \Delta u_1  +  \epsilon (t) (u_1 -u_2) = h_1 \ \ on  \ \Omega \times ]0, +\infty[\\
\rule{0pt}{25pt}
\frac{\partial^2 u_2}{\partial t^2} + \gamma \frac{\partial u_2}{\partial t}  - \alpha_2 \Delta u_2  +  \epsilon (t) (u_2 -u_1) = h_2 \  \ on  \ \Omega \times ]0, +\infty[\\
\rule{0pt}{25pt}
\frac{\partial u_1}{\partial n}= \frac{\partial u_2}{\partial n}=0  \quad on \ \partial \Omega \times ]0, +\infty[\\
\rule{0pt}{25pt}
u_1 (0) =u_{01}, u_2(0)= u_{02}\\
\rule{0pt}{25pt}
\frac{\partial u_1}{\partial t}(0)= v_{01}, \ \frac{\partial u_2}{\partial t}(0)= v_{02}.
\end{array}\right.
\end{equation}

\smallskip

Let us interpret the system (\ref{couplwave1}) in a functional setting. Take

\begin{itemize}
\item $\mathcal H = L^2 (\Omega) \times L^2 (\Omega) $ with its Hilbert product structure.

\smallskip
 
\item  $\Phi:\mathcal H = L^2 (\Omega) \times L^2 (\Omega) \to \mathbb R \cup \{+\infty\}$, with for any $v=(v_1,v_2)  \in \mathcal H$, \
$\Phi(v) = \Phi_1(v_1) + \Phi_2(v_2) $, 
\begin{equation}\label{couplwave2}
\Phi_i (v) =\left\{
\begin{array}{l}
\frac{\alpha_i}{2}   \int_{\Omega} |\nabla v(x)|^2 dx -  \int_{\Omega} h_i (x) v(x) dx    \quad if \    v\in  H^1 (\Omega) \\
\rule{0pt}{25pt}
+ \infty \quad if \ v\in L^2 (\Omega) \setminus H^1 (\Omega)
\end{array}\right.
\end{equation}
\item : $\Psi: \mathcal H = L^2 (\Omega) \times L^2 (\Omega) \to \mathbb R$, with for any $v=(v_1,v_2)  \in \mathcal H$,
$$\Psi (v)= \frac{1}{2}\int_{\Omega} |v_1(x)-v_2 (x)|^2 dx. $$
\end{itemize}

\smallskip

The system (\ref{couplwave1}) can be equivalently formulated as the Cauchy  problem for the following evolution equation
$$
   \ddot{u}(t)   + \gamma  \dot{u}(t) +    \partial\Phi (u(t)) +   \epsilon (t) \nabla \Psi  (u(t)) = 0. 
$$
In the above equation, $\partial \Phi$ is the subdifferential of $\Phi$ in $\mathcal H$, which, because of the decoupled formulation of $\Phi$, gives for each component, the (minus) Laplace operator, with homogenous Neumann boundary condition, and with respective coefficients
 $\alpha_1$ and $\alpha_2$,  see 
for example \cite[Theorem 17.2.12]{ABM}. Note that computing the gradient at $u$ of the coupling potential $\Psi $ gives two opposite forces $u_1 -u_2$ and $u_2 -u_1 $ (action and reaction principle).

As we have already noted, Theorem \ref{basicthm}  can not be  directly  applied to this sytem, because $\Phi$ is only lower semicontinuous. Still, our approach can be adapted to this situation.
The point is that the dimension of the space is not involved in our analysis, and  Rellich-Kondrakov compact embedding of
$H^1 (\Omega)$ into $L^2 (\Omega) $, makes the functional $\Phi$  inf-compact on $\mathcal H = L^2 (\Omega) \times L^2 (\Omega) $.
Thus we are naturally led to develop a similar analysis to that used in Theorem \ref{basicthm}.
We omit the details of the proof of convergence, which combines our proof  with  the arguments used in the study of the wave equation in \cite[Theorem 2.1]{AA}, and \cite[Proposition 4.1]{AtCz1}. 
Let us just describe the asymptotic selection result.
 
Because of the compatibility condition $\int_{\Omega}h_i=0$, the  variational problem

\begin{equation}\label{Neumann1}
 \min \left\lbrace \frac{\alpha_i}{2}   \int_{\Omega} |\nabla v(x)|^2 dx -  \int_{\Omega} h_i (x) v(x) dx :    \quad   v\in  H^1 (\Omega), \ \int_\Omega v(x)dx =0 \right\rbrace 
\end{equation}
has a unique solution, which is denoted by $\bar{u}_i$, see for example \cite[Theorem 6.2.3]{ABM}. As a consequence, 
the solution set $C= \argmin \Phi$ is given by
\begin{equation}\label{Neumann2}
 C = \left\lbrace  (\bar{u}_1 + r_1, \bar{u}_2 + r_2): \ r_1\in \mathbb R, \ \ r_2\in \mathbb R \right\rbrace  .
\end{equation}
 Because of the quadratic property of $\Psi$, the moderate decrease condition 
 $({\mathcal H}_1)$ can be equivalently written 
\begin{equation}\label{Neumann3}
 \epsilon (\cdot) \in L^2 (0, + \infty) \setminus L^1 (0, + \infty)..
\end{equation}
Thus, under the above condition, we obtain that the solution  trajectory $t \mapsto (u_1 (t), u_2 (t))$ of (\ref{couplwave1}) converges strongly in $\mathcal H = L^2 (\Omega) \times L^2 (\Omega)$ to some $u_{\infty}= (u_{1,\infty}, u_{2,\infty})    $   which is a solution of the minimization problem
\begin{equation}\label{Neumann4}
\min \left\lbrace  \int_{\Omega} |v_1(x)-v_2 (x)|^2 dx:
 \  v_i= \bar{u}_i + r_i , r_i \in \mathbb R \right\rbrace .
\end{equation}
An elementary calculus gives $u_{1,\infty}= \bar{u}_1 + r_1    $,
$u_{2,\infty}= \bar{u}_2 + r_2    $, with $r_1=r_2$. Equivalently,
$$
\int_{\Omega} u_{1,\infty} (x) dx = \int_{\Omega} u_{2,\infty} (x)dx 
$$
which means that, at the limit, the two waves stabilize at equilibria that are as close as possible, in the sense that they have the same mean value.
Convergence of the energies allows to pass from the convergence for the $L^2$ norm to the convergence for the $H^1$ norm.

This model suggests new developments. For example, it would be interesting to study the case where the two waves propagate in different domains, and the coupling occurs only on their common boundary.

\end{document}